\documentclass[11pt]{article}

\usepackage{graphicx,tikz,color,url}
\usepackage{amsmath,amssymb,latexsym,amsfonts}
\usepackage{mathtools}
\usepackage{xspace}
\usepackage[T1]{fontenc}
\usepackage{cite}

\newtheorem{theorem}{Theorem}[section]

\newtheorem{lemma}[theorem]{Lemma}
\newtheorem{proposition}[theorem]{Proposition}

\newtheorem{corollary}[theorem]{Corollary}

\newcommand{\proof}{\noindent{\bf Proof.\ }}
\newcommand{\qed}{\hfill $\square$ \bigskip}
\newcommand{\cp}{\,\square\,}

\newcommand{\diam}{{\rm diam}}
\newcommand{\gp}{{\rm gp}}
\newcommand{\sr}{{}_{\rm SR}}

\newcommand{\gpt}{\gp_{\rm t}}
\newcommand{\gpd}{\gp_{\rm d}}
\newcommand{\gpo}{\gp_{\rm o}}

\newcommand{\name}{$P_4$-inner isometric}

\title{Variety of general position problems in graphs}

\author{
Jing Tian $^{a,b}$ \and Sandi Klav\v{z}ar $^{b,c,d}$\\\\
$^{a}$ \small School of Science, Zhejiang University of Science and Technology, \\
\small Hangzhou, Zhejiang 310023, PR China\\
\small {\tt jingtian526@126.com}\\
$^{b}$ \small Institute of Mathematics, Physics and Mechanics, Ljubljana, Slovenia \\
$^{c}$\small Faculty of Mathematics and Physics, University of Ljubljana, Slovenia\\
$^{d}$ \small Faculty of Natural Sciences and Mathematics, University of Maribor, Slovenia\\
\small {\tt sandi.klavzar@fmf.uni-lj.si}\\
}
\date{}

\begin{document}
\maketitle

\begin{abstract}
Let $X$ be a vertex subset of a graph $G$. Then $u, v\in V(G)$ are $X$-positionable if $V(P)\cap X \subseteq \{u,v\}$ holds for any shortest $u,v$-path $P$. If each two vertices from $X$  are $X$-positionable, then $X$ is a general position set. The general position number of $G$ is the cardinality of a largest general position set of $G$ and has been already well investigated. In this paper a variety of general position problems is introduced based on which natural pairs of vertices are required to be $X$-positionable. This yields the total (resp.\ dual, outer) general position number. It is proved that the total general position sets coincide with sets of simplicial vertices, and that the outer general position sets coincide with sets of mutually maximally distant vertices. It is shown that a general position set is a dual general position set if and only if its complement is convex. Several sufficient conditions are presented that guarantee that a given graph has no dual general position set.  The total general position number, the outer general position number, and the dual general position number of arbitrary Cartesian products are determined. 
\end{abstract}


\noindent
AMS Subj.\ Class.\ (2020): 05C12, 05C69, 05C76

\maketitle

\section{Introduction}

General position sets were introduced to graph theory independently in~\cite{Ullas-2016, Manuel-2018}, but in the special case of hypercubes, these sets had been studied previously in~\cite{Korner-1995}. Finding a largest general position set in a graph $G$ is an NP-hard problem~\cite{Manuel-2018}. In~\cite{Anand-2019}, a characterization of general position sets in graphs was proved and the general position number for the complement of bipartite graphs, the complement of trees, and the complement of hypercubes was determined. Afterwards, general position sets received a wide attention, see~\cite{Ghorbani-2021, Irsic-2023, Klavzar-Rall-Yero-2021, Klavzar-Tian-2023, Manuel-2022} for a selection of recent developments.

Mutual-visibility sets are closely related to general position sets and were recently introduced by Di Stefano~\cite{DiStefano-2022}. While investigating mutual-visibility sets in strong products, there was a natural need to introduce total mutual-visibility sets~\cite{cicerone-2023+}. This has led to the variety of mutual-visibility sets which was formalized in~\cite{CiDiDrHeKlYe-2023}. Motivated by the variety of mutual-visibility sets,  the main purpose of this paper is to introduce and study the corresponding variety of general position sets. It consists of general position sets, total general position sets, outer general position sets, and dual general position sets.

The rest of the paper is organized as follows.  In the rest of the introduction we give definitions needed. In the next section we introduce the variety of the general position sets and characterize total general position sets and outer general position sets. In Section~\ref{sec:dual} we focus on dual general position sets. We first observe that a general position set is a dual general position set if and only if its complement is convex. Then we consider graphs with small dual general position numbers. In particular, several sufficient conditions are presented that guarantee that a given graph has no dual general position set at all. In Section~\ref{sec:in-Cartesian} we determine the total general position number, the outer general position number, and the dual general position number of arbitrary Cartesian products. 

For a natural number $n$ we set $[n] = \{1,\ldots, n\}$. All graphs $G = (V(G), E(G))$ in the paper are connected unless otherwise stated.  The {\em order} of $G$ is the value of $|V(G)|$. The open neighborhood of a vertex $u$ of $G$ is denoted by $N_G(u)$.  The {\em degree} of a vertex $u$ of $G$ is $\deg_G(u) = |N_G(u)|$. If $X\subseteq V(G)$, then $G[X]$ denotes the subgraph of $G$ induced by $X$. Moreover, $G-X$ is the subgraph of $G$ obtained from $G$ by deleting all vertices from $X$. If $G$ is not a tree, then its {\em girth} $g(G)$ is the length of a shortest cycle of $G$. A vertex of $G$ is {\em simplicial} if its neighbourhood induces a complete subgraph. The set of simplicial vertices of $G$ will be denoted by $S(G)$ and the cardinality of $S(G)$ by $s(G)$. The clique number of $G$ is denoted by $\omega(G)$. 

The {\em distance} $d_G(u,v)$ between vertices $u$ and $v$ of $G$ is the usual shortest-path distance. The {\em diameter} $\diam(G)$ of $G$ is the maximum distance between pairs of vertices of $G$. A subgraph $G'$ of a graph $G$ is \emph{isometric}, if for every two vertices $x$ and $y$ of $G'$ we have $d_{G'}(x,y) = d_{G}(x,y)$. A subgraph $G'$ of a graph $G$ is  \emph{convex}, if for every two vertices of $G'$, every shortest path in $G$ between them lies completely in $G'$. By abuse of language we also say that a set of vertices is convex if it induces a convex subgraph.

The {\em Cartesian product} $G\cp H$ of graphs $G$ and $H$ has the vertex
set $V(G\cp H) = V(G)\times V(H)$, and vertices $(g, h$) and $(g', h')$ are adjacent if either $gg'\in E(G)$ and $h = h'$, or $g = g'$ and $hh'\in E(H)$. The {\em direct product} $G\times H$ has the same vertex set as $G\cp H$, while vertices $(g, h$) and $(g', h')$ are adjacent if $gg'\in E(G)$ and $hh'\in E(H)$. The {\em strong product} $G\boxtimes H$ also has the same vertex set as $G\cp H$, while $E(G\boxtimes H) = E(G\cp H) \cup E(G\times H)$. The {\em join} $G\oplus H$ is the graph obtained from the disjoint union of $G$ and $H$ by adding all possible edges between vertices from $G$ and vertices from $H$.

\section{The variety}
\label{sec:variety}

In this section we first introduce the announced variety of general position sets. After that we characterize total general position sets and outer general position sets.

Let $G = (V(G), E(G))$ be a graph and $X\subseteq V(G)$. Vertices $u,v\in V(G)$ are {\em $X$-positionable} if for any shortest $u,v$-path $P$ we have $V(P)\cap X \subseteq \{u,v\}$. Note that each pair of adjacent vertices is $X$-positionable. Set $\overline{X} = V(G)\setminus X$. Then we say that $X$ is 
\begin{itemize}
\item a \emph{general position set}, if every $u,v\in X$ are $X$-positionable;
\item a \emph{total general position set}, if every $u,v\in V(G)$ are $X$-positionable;
\item an \emph{outer general position set}, if every $u,v\in X$ are $X$-positionable, and every $u\in X$, $v\in \overline{X}$ are $X$-positionable; and
\item a \emph{dual general position set}, if every $u,v\in X$ are $X$-positionable, and every $u,v\in \overline{X}$ are $X$-positionable.
\end{itemize}
The cardinality of a largest general position set, a largest total general position set, a largest outer general position set, and a largest dual general position set will be respectively denoted by $\gp(G)$, $\gpt(G)$, $\gpo(G)$, and $\gpd(G)$. Also, these graph invariants will be respectively called the {\em general position number}, the {\em total general position number}, the {\em outer general position number}, and the {\em dual general position number} of $G$.
Moreover, for any invariant $\tau(G)$ from the above ones, by a {\em $\tau$-set} we mean any set of vertices of cardinality $\tau(G)$. In addition, for any two invariants $\tau_1(G)$ and $\tau_2(G)$, by a {\em $(\tau_1,\tau_2)$-graph} we mean any graph $G$ with $\tau_1(G)=\tau_2(G)$.

If $G$ is a graph, then by definition,
\begin{align}
\gp(G) & \ge \gpo(G) \ge \gpt(G)\quad {\rm and} \label{eq:1}\\
\gp(G) & \ge \gpd(G) \ge \gpt(G)\label{eq:2}\,.
\end{align}

If $G$ is a block graph, then $\gp(G) = s(G)$, see~\cite[Theorem 3.6]{Manuel-2018}. Moreover, we can check directly that the set of simplicial vertices of $G$ is also a total general position set. Hence block graphs are $(\gp,\gpt)$-graphs. Having~\eqref{eq:1} and~\eqref{eq:2} in mind we thus get that if $G$ is a block graph, then
$$\gp(G) = \gpo(G) = \gpd(G) = \gpt(G) = s(G)\,.$$
In particular, if $n\ge 2$, then $\tau(P_n) = 2$ for each $\tau\in \{\gp, \gpo, \gpd, \gpt\}$. Moreover,
\begin{itemize}
\item
 each pair of distinct vertices of $P_n$ forms a $\gp$-set of $P_n$,
\item
 $\{1,n\}$ is the only $\gpo$-set of $P_n$,
\item
 $\{1,2\}$, $\{1,n\}$, and $\{n-1,n\}$ are the only $\gpd$-sets of $P_n$,
\item
 $\{1,n\}$ is the only $\gpt$-set of $P_n$.
\end{itemize}

Contrary to block graphs, by considering Cartesian products of two complete graphs we will see that all these four parameters can be pairwise arbitrary different.

We now characterize total general position sets as follows.

\begin{theorem}
\label{thm:total-characterization}
Let $G$ be a graph and $X\subseteq V(G)$. Then $X$ is a total general position set of $G$ if and only if $X\subseteq S(G)$. Moreover, $\gpt(G)=s(G)$.
\end{theorem}

\proof
Assume first that $X$ is a total general position set of $G$ and consider an arbitrary vertex $x$ from $X$.
Suppose that $x\not\in S(G)$. Then $x$ has two adjacent vertices $y$ and $z$ such that $d_G(y,z)=2$, which means that there exists a shortest $y,z$-path containing $x$. Hence the vertices $y$ and $z$ are not $X$-positionable, a contradiction. As a consequence, we conclude that $x$ must belong to $S(G)$ and then $X\subseteq S(G)$.

To prove the converse, assume that $X\subseteq S(G)$. Consider any two vertices $u$ and $v$ from $V(G)$ and let $P_{uv}$ be a shortest $u,v$-path of $G$. If $u$ is adjacent to $v$, then $V(P_{uv}) = \{u,v\}$ and hence $u$ and $v$ are $X$-positionable. Assume next that $d_G(u,v)\geq 2$. Let $P_{uv}$ be the path $u=x_1, x_2, \ldots, x_k=v$. Then $k\geq 3$.
Suppose that $x_i\in X$, where $2\leq i\leq k-1$. Since $x_i$ is a simplicial vertex of $G$, $x_{i-1}$ must be adjacent to $x_{i+1}$ and hence $x_1,\ldots, x_{i-1},x_{i+1},\ldots, x_k$ is a $u,v$-path shorter than $P_{uv}$. Since this is not possible, $V(P_{uv})\cap X\subseteq \{u,v\}$, and thus $u$ and $v$ are $X$-positionable. Hence $X$ is a total general position set of $G$, and then we can conclude that $\gpt(G)=s(G)$.
\qed

\begin{corollary}
\label{cor:gpt=0}
If $G$ is a graph, then $\gpt(G)=0$ if and only $G$ has no simplicial vertices.
\end{corollary}

We continue with a characterization of outer general position sets. For this sake, the following concepts are crucial. A vertex $u$ of a graph $G$ is \emph{maximally distant} from $v\in V(G)$ if for every $w\in N_G(u)$ it holds $d_G(v,w)\le d_G(u,v)$. If also $v$ is maximally distant from $u$, then $u$ and $v$ are said to be \emph{mutually maximally distant}. The {\em strong resolving graph} $G\sr$ of $G$ has $V(G\sr) = V(G)$ and two vertices being adjacent in $G\sr$ if they are mutually maximally distant in $G$. This notion was introduced in~\cite{Oellermann-2007} as a tool to study the strong metric dimension. The paper~\cite{Kuziak-2018} gives a survey on strong resolving graphs with an emphasize on the realization problem (which graphs have a given graph as their strong resolving graph) and the characterization problem (characterize graphs that are strong resolving graphs of some graphs). From our perspective, in~\cite{Klavzar-Yero-2019} it was proved that $\gp(G)\ge \omega(G\sr)$ holds for any connected graph $G$.

Now we can characterize outer general position sets as follows.

\begin{theorem}
\label{thm:characterize outer}
Let $G$ be a connected graph and $X\subseteq V(G)$, $|X|\geq 2$. Then $X$ is an outer general position set of $G$ if and only if each pair of vertices from $X$ is mutually maximally distant. Moreover,
$$\gpo(G) = \omega(G\sr)\,.$$
\end{theorem}

\proof
Let $X\subseteq V(G)$ be a set with $|X|\ge 2$.

Assume first that $X$ is an outer general position set of $G$ and let $x,y\in X$. If $x$ and $y$ are not mutually maximally distant, there exists (at least) one neighbor of $x$ or $y$, say $w\in N_G(x)$, such that $d_G(w,y) = d_G(x,y) + 1$. Hence there exists a shortest $w,y$-path which contains $x$. Then no matter whether $w\in X$ or $w\notin X$, this contradicts the fact that $\{x,y\}\subseteq X$ is an outer general position set.

Conversely, assume that any two vertices from $X$ are mutually maximally distant. Consider any two vertices $u$ and $v$. If $u,v\in X$, the internal vertices of all shortest $u,v$-paths are not from $S$ as $u$ and $v$ are mutually maximally distant. Hence $u$ and $v$ are $X$-positionable. Assume next that, without loss of generality, $u\in X$ and $v\in V(G)\setminus X$. If there exists a vertex $w\in X\setminus \{u\}$ lying on a shortest $u,v$-path, then $u$ and $w$ are not mutually maximally distant. Hence $u,v$ are $X$-positionable.

The formula $\gpo(G) = \omega(G\sr)$ now follows by the above characterization of outer general position sets and by the definition of the strong resolving graph.
\qed

\begin{corollary}
If $G$ is a connected graph of order at least $2$, then $\gpo(G)\ge 2$.
\end{corollary}

\proof
Let $u$ and $v$ be vertices of $G$ such that $d_G(u,v) = \diam(G)$. Then $u$ and $v$ are mutually maximally distant, hence $\gpo(G)\ge 2$ by Theorem~\ref{thm:characterize outer}.
\qed

As already mentioned, it was proved in~\cite{Klavzar-Yero-2019} that $\gp(G)\ge \omega(G\sr)$ holds for any connected graph $G$. Hence by Theorem~\ref{thm:characterize outer} we have
$$\gp(G)\ge \omega(G\sr) = \gpo(G)\,.$$
Thus, knowing that $\gp(G) = \omega(G\sr)$ holds for some graph $G$, we also know $\gpo(G)$. For instance, \cite[Proposition 4.5]{Klavzar-Yero-2019} asserts that if $r_1\ge t_1\ge 1$ and $r_2\ge t_2\ge 1$, 
then $$\gp(K_{r_1,t_1}\boxtimes K_{r_2,t_2}) = r_1r_2 = \omega((K_{r_1,t_1}\boxtimes K_{r_2,t_2})\sr)\,,$$
which in turn implies that
$$\gpo(K_{r_1,t_1}\boxtimes K_{r_2,t_2}) = r_1r_2\,.$$

\section{Dual general position sets}
\label{sec:dual}

From Theorem~\ref{thm:total-characterization} it follows that if $X$ is a total general position set of $G$, then each subset of $X$ is also a total general position set of $G$. Similarly, by Theorem~\ref{thm:characterize outer} this hereditary property also holds for outer general position sets. In addition, this property is also known to hold for general position sets. A bit surprisingly,  if $X$ is a dual general position set of $G$ and $Y\subseteq X$, then $Y$ need not be a dual general position set of $G$. To see this, consider two adjacent vertices of $C_5$ which form a dual general position set, however, one vertex of $C_5$ does not form such a set. So each of the properties of being in general position, being in total general position, and being in outer general position is hereditary, but the property of being in dual general position property is not.

The above consideration indicates that the dual general position is intrinsically different from the other three invariants. In this section we have a closer look to it. We first characterize which general position sets are dual general position sets. Then we respectively  consider graphs $G$ with $\gpd(G) = 0$, $\gpd(G) = 1$, and $\gpd(G) \ge 2$.

Let $X$ be a dual general position set in a graph $G$. Then, clearly, $X$ is a general position set. Hence to find a largest dual general position set in $G$ if suffices to check all general position sets in $G$ and check if they are also dual. For this task, the following result is useful.

\begin{theorem}
\label{thm:dual-characterization}
Let $X$ be a general position set of a graph $G$. Then $X$ is a dual general position set if and only if $G-X$ is convex.
\end{theorem}

\proof
Assume that $X$ is a dual general position set. Let $x,y\in V(G)\setminus X$  and let $P$ be a shortest $x,y$-path. Since $X$ is a dual general position set, the vertices $x$ and $y$ are $X$-positionable which in turn implies that $V(P)\cap X = \emptyset$. It follows that $G-X$ is convex.

Conversely, assume that $X$ is a general position set and $G-X$ is convex. If $x,y\in X$, then $x$ and $y$ are $X$-positionable since $X$ is a general position set. Consider next $x,y\in V(G)\setminus X$. Since $G-X$ is convex, each shortest $x,y$-path lies in $G-X$ hence no such path contains vertices from $X$. We can conclude that $X$ is a dual general position set.
\qed

For the later use we state explicitly the following consequence of Theorem~\ref{thm:dual-characterization}.

\begin{corollary}
\label{cor:simplicial-are-dual}
If $G$ is a graph and $X\subseteq S(G)$, then $X$ is a dual general position set.
\end{corollary}

\proof
Since $X\subseteq S(G)$ we see that $X$ is a general position set and also that $G-X$ is convex. Hence the assertion follows by Theorem~\ref{thm:dual-characterization}.
\qed

\subsection{Graphs $G$ with $\gpd(G) = 0$}
\label{sec:gpd=0}

In this subsection we focus on graphs $G$ with $\gpd(G) = 0$. It follows from Theorem~\ref{thm:dual-characterization} that $\gpd(G) = 0$ if and only if for every general position set $X$ the subgraph  $G-X$ is not convex in $G$. Since dual general position sets are not hereditary (recall the example of $C_5$), we must consider all general position sets, not only singletons.

We say that an edge $e$ of a graph $G$ is {\em \name}, if $e$ is the middle edge of some isometric $P_4$. 

\begin{proposition}
\label{pro:gp_d=0}
Let $G$ be a connected graph. If each edge of $G$ is \name, then $\gpd(G)=0$.
\end{proposition}

\proof
Assume that $G$ is a graph in which each edge is \name. Then clearly $\delta(G)\geq 2$. Suppose on the contrary that $\gpd(G)\geq 1$ and let $x$ be a vertex of $G$ from some dual general position set $X$.
Let $y$ be a neighbor of $x$. From our assumption, the edge $xy$ is \name\ and let $x', x, y, y'$ be the vertices of an isometric $P_4$, denote it by $Q$, where $x'x\in E(G)$ and $yy'\in E(G)$. Then at least one of the vertices $x'$ and $y$ belongs to $X$, for otherwise $x'$ and $y$ are not $X$-positionable.

Suppose first that $y\in X$. It implies that $y'\not\in X$ for otherwise $x$ and $y'$ are not $X$-positionable. Analogously, $x'\not\in X$ also holds. But then the vertices $x'$ and $y'$ are not $X$-positionable.

Suppose second that $x'\in X$. Then $y\notin X$, for otherwise $y$ and $x'$ are not $X$-positionable. Consider now an isometric $P_4$, say $R$, such that the edge $x'x$ is the middle edge of $R$. Let $R$ be the path on the vertices $z', x', x, z$, where $z'x'\in E(G)$ and $xz\in E(G)$. Since $x,x'\in X$ we see that $z'\notin X$ and $z\notin X$. Since $Q$ is isometric, we have $z'\ne y$ and $z'\ne y'$. If $z=y$, then $z'$ and $z=y$ are not $X$-positionable.  Similarly, if $ z\ne y$, then again $z'$ and $z$ are not $X$-positionable. This final contradiction implies that $\gpd(G) = 0$.
\qed

\begin{corollary}
\label{cor:gpd=0 of girth 6}
Let $G$ be a connected graph with $g(G) \ge 6$. Then $\gpd(G)=0$ if and only if $\delta(G)\geq 2$.
\end{corollary}

\proof
If $\delta(G) = 1$ and $u$ is a pendant vertex of $G$, then $\{ u \}$ is a dual general position set, hence $\gpd(G)\ge 1$.

To prove the other direction, assume that $\delta(G)\geq 2$. Since $g(G)\ge 6$ we see that each edge of $G$ is \name. Proposition~\ref{pro:gp_d=0} yields the conclusion.
\qed

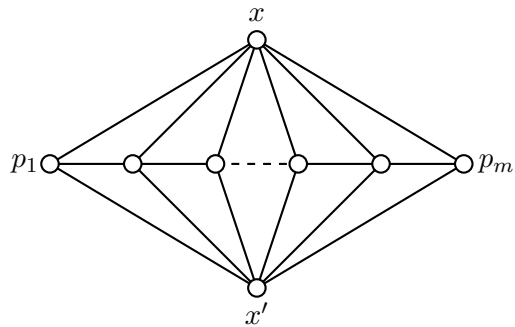
\begin{figure}[ht!]
\begin{center}
\begin{tikzpicture}[scale=1.1,style=thick]
\tikzstyle{every node}=[draw=none,fill=none]
\def\vr{3pt}

\begin{scope}[yshift = 0cm, xshift = 0cm]
    \node [below=0.5mm] at (0,1.5) {};
    \node [below=0.5mm] at (0,-1.5){};
    \node [below=0.5mm] at (0.5,0) {};
    \node [below=0.5mm] at (1.5,0) {};
    \node [below=0.5mm] at (2.5,0) {};
    \node [below=0.5mm] at (-0.5,0) {};
    \node [below=0.5mm] at (-1.5,0) {};
    \node [below=0.5mm] at (-2.5,0) {};

\path (-2.5,0) coordinate (x1);
\path (-1.5,0) coordinate (x2);
\path (-0.5,0) coordinate (x3);
\path (0.5,0) coordinate (x4);
\path (1.5,0) coordinate (x5);
\path (2.5,0) coordinate (x6);
\path (0,1.5) coordinate (x7);
\path (0,-1.5) coordinate (x8);

\draw (x1)--(x2)--(x3);
\draw (x4)--(x5)--(x6);
\draw (x7)--(x1)--(x8);
\draw (x7)--(x2)--(x8);
\draw (x7)--(x3)--(x8);
\draw (x7)--(x4)--(x8);
\draw (x7)--(x5)--(x8);
\draw (x7)--(x6)--(x8);
\draw [dashed](x3)--(x4);

\draw (x1)  [fill=white] circle (\vr);
\draw (x2)  [fill=white] circle (\vr);
\draw (x3)  [fill=white] circle (\vr);
\draw (x4)  [fill=white] circle (\vr);
\draw (x5)  [fill=white] circle (\vr);
\draw (x6)  [fill=white] circle (\vr);
\draw (x7)  [fill=white] circle (\vr);
\draw (x8)  [fill=white] circle (\vr);

\draw (0,1.8) node {$x$};
\draw (0,-1.8) node {$x'$};
\draw (-2.8,0) node {$p_1$};
\draw (2.9,0) node {$p_m$};

\end{scope}
\end{tikzpicture}
\end{center}
\caption{The graph $G_m$.}
\label{fig: graph 1}
\end{figure}

Next, we give an infinite family of graphs whose dual general position number is zero, yet none of their edges is \name. Let $m\ge 5$ and consider the join graph $G_m = P_m \oplus 2K_1$, see Fig.~\ref{fig: graph 1}.

It is readily verified that no edge of $G_m$ is \name. On the other hand, the following holds.

\begin{proposition}
\label{pro:dual-paths join two isolated vertices}
If $m\geq 5$, then $\gpd(G_m) = 0$.
\end{proposition}

\proof
Set $G_m=P_m\oplus 2K_1$, let $V(P_m)=\{p_1,\ldots,p_m\}$ with natural adjacencies, and let $V(2K_1)=\{x,x'\}$, see Fig.~\ref{fig: graph 1} again.

Let $X$ be an arbitrary dual general position set of $G_m$. If $x, x'\not\in X$, then $X = \emptyset$ for otherwise $x$ and $x'$ are not $X$-positionable. In the rest we may hence assume without loss of generality that $x\in X$, for otherwise we are done.

We first claim that $x'\notin X$. Indeed, otherwise no vertex $p_i$, $i\in [m]$, lies in $X$, but then no two vertices $p_i$ and $p_j$, $|i - j|\ge 2$, are not $X$-positionable.  Further and similarly, we have $\{p_1,\ldots,p_m\}\cap X\neq \emptyset$. Hence there exists $i\in [m]$ such that $p_i\in X$. Select and fix $i$ to be the smallest such index. Then by the symmetry we may assume that $1\le i\le \lceil m/2\rceil$. We now distinguish two cases.

Suppose first that $i=1$. Then it follows that $p_j\not\in X$ for $3\leq j\leq m$. Since $m\geq 5$, the vertex $x$ lies in the middle of a shortest $p_3,p_m$-path, hence the vertices $p_3$ and $p_m$ are not $X$-positionable.

Suppose second that $2\leq i\le \lceil m/2\rceil$. By the way $i$ is selected and by the symmetry we have $p_1\notin X$ and $p_{m}\notin X$. But then $p_1$ and $p_m$ are not $X$-positionable.
\qed

For the cycles $C_4$ and $C_5$, it is observed that $\gpd(C_4)=\gpd(C_5)=2$. These graphs can be considered as special cases of generalized theta graphs which are defined as follows. For positive integers $1\leq \ell_1\leq \cdots\leq \ell_k$ and $\ell_2\geq 2$,  the \textit{generalized theta graph} $\Theta(\ell_1,\ldots,\ell_k)$ is obtained by joining two vertices $a$ and $b$ with $k$ internally disjoint paths of lengths $\ell_1,\ldots, \ell_k$.

\begin{proposition}
\label{pro:dual for theta graph}
Let $1\leq \ell_1\leq \cdots\leq \ell_k$, where $k\geq 2$, $\ell_2\geq 2$. Then $\gpd(\Theta(\ell_1,\ldots,\ell_k)) = 0$ if and only if one of the following cases holds:
\begin{enumerate}
\item[(i)] $k=2$, $\ell_1+\ell_2\geq 6$;
\item[(ii)] $k\geq 3$, $\ell_1=1$, $\ell_2\geq 5$;
\item[(iii)] $k\geq 3$, $\ell_1=2$, and $\ell_i\neq 3$ for $2\leq i\leq k$; \item[(iv)] $k\geq 3$, $\ell_1\geq 3$.
\end{enumerate}
\end{proposition}

\proof
Set $\Theta=\Theta(\ell_1,\ldots,\ell_k)$ for the rest of the proof and let $Q_1,\ldots, Q_k$ be the internally disjoint paths of lengths $\ell_1,\ldots, \ell_k$ of $\Theta$ connecting $a$ and $b$.

$(i)$ If $k=2$ and $\ell_1+\ell_2\leq 5$, then $\Theta\in \{C_3, C_4, C_5\}$, hence $\gpd(\Theta)\geq 2$. On the other hand, if $\ell_1+\ell_2\geq 6$, then $\Theta$ is a cycle of order at least $6$ and thus each edge of $\Theta$ is \name. Hence $\gpd(\Theta) = 0$ by Proposition~\ref{pro:gp_d=0}.

$(ii)$ Let $k\geq 3$ and $\ell_1=1$. Assume first that $\ell_2\leq 4$. If $\ell_2 = 2$, then the middle vertex of $Q_2$ is simplicial and by Corollary~\ref{cor:simplicial-are-dual}, $\{x\}$ is a dual general position set and so $\gpd(\Theta)\geq 1$. If $\ell_2 = 3$, then $\Theta[Q_1\cup Q_2]\cong C_4$. Let $u$ and $v$ be the internal vertices of $Q_2$. Then we infer that $\Theta-\{u,v\}$ is convex and using Theorem~\ref{thm:dual-characterization} we get $\gpd(\Theta)\geq 2$. Similarly, we get that $\gpd(\Theta)\geq 2$ if $\ell_2 = 4$. Assume finally that $\ell_2\geq 5$. Then each edge of $\Theta$ is \name. By Proposition~\ref{pro:gp_d=0} we thus have $\gpd(\Theta) = 0$.

$(iii)$ Let $k\geq 3$ and $\ell_1=2$. Assume first that there exists an index $i$ such that $\ell_i = 3$ and let $i$ be the smallest such index. Note that $i\ge 2$. If $x_i$ and $x_i'$ are the two internal vertices of $Q_i$, then $\Theta - \{x_i,x_i'\}$ is convex, hence Theorem~\ref{thm:dual-characterization} implies that $\gpd(\Theta)\geq 2$.

Assume second that $\ell_i \neq 3$ for each $2\leq i\leq k$. If $\ell_1 = \cdots = \ell_k = 2$, then $\Theta \cong K_{2,k}$. A possible dual general position set wound need to contain two adjacent vertices of each $4$-cycle, but this is not possible. So $\gpd(K_{2,k}) = 0$ for $k\ge 3$. Let next $j$ be the smallest index such that $\ell_j \ge 4$, so that $\ell_1 = \cdots = \ell_{j-1} = 2$. (It is possible that $j=2$.) Then the subgraph of $\Theta$ induced by $Q_1\cup \cdots \cup Q_{j-1}$ is isomorphic to $K_{2,j-1}$ and we readily see that no vertex from it can lie in a dual general position set. In addition, the vertices from $Q_1\cup Q_{j'}$, where $j'\ge j$, induce an isometric cycle of $\Theta$ of order at least $6$, from which we conclude that none of the vertices from the cycle can lie in a dual general position set. We conclude that $\gpd(\Theta) = 0$.

$(iv)$ In this case we have $g(\Theta) \ge 6$, hence the assertion follows by Corollary~\ref{cor:gpd=0 of girth 6}.
\qed

\subsection{Graphs $G$ with $\gpd(G) = 1$}
\label{sec:gpd=1}

We next consider graphs $G$ with $\gpd(G) = 1$.

\begin{proposition}
\label{prop:dual=1}
If $G$ is a connected graph with $\gpd(G) = 1$, then $s(G) = 1$.
\end{proposition}

\proof
Let $G$ be a connected graph with $\gpd(G) = 1$ and let $\{x\}$ be a $\gpd$-set of $G$. Suppose first that $s(G) = 0$. Then it follows that the order of $G$ at least $4$ as each of the smaller graphs contains at least one simplicial vertex. In particular, since $x$ is not simplicial, it has two neighbors, say $x'$ and $x''$, such that $d_G(x',x'')=2$. But then $x'$ and $x''$ are not $\{x\}$-positionable. This contradiction implies that $s(G)\ge 1$. On the other hand, Corollary~\ref{cor:simplicial-are-dual} yields that $s(G)\le 1$ and we are done.
\qed

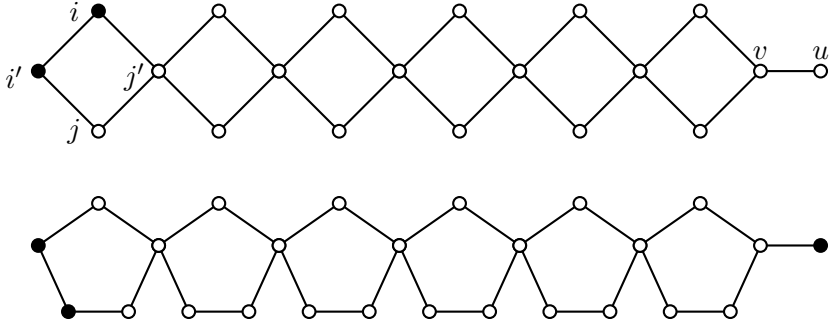
\begin{figure}[ht!]
\begin{center}
\begin{tikzpicture}[scale=0.8,style=thick,x=1cm,y=1cm]
\def\vr{3pt}

\begin{scope}[xshift=0cm, yshift=0cm] 
\coordinate(x1) at (1,3);
\coordinate(x2) at (2,4);
\coordinate(x3) at (1,5);
\coordinate(x4) at (0,4);
\draw (x1) -- (x2) -- (x3) -- (x4) --(x1);
\draw(x1)[fill=white] circle(\vr);
\draw(x2)[fill=white] circle(\vr);
\draw(x3)[fill=black] circle(\vr);
\draw(x4)[fill=black] circle(\vr);
\node at (0.6,3) {$j$};
\node at (0.6,5) {$i$};
\node at (-0.4,3.9) {$i'$};
\node at (1.6,3.9) {$j'$};
\end{scope}

\begin{scope}[xshift=2cm, yshift=0cm] 
\coordinate(x1) at (1,3);
\coordinate(x2) at (2,4);
\coordinate(x3) at (1,5);
\coordinate(x4) at (0,4);
\draw (x1) -- (x2) -- (x3) -- (x4) --(x1);
\draw(x1)[fill=white] circle(\vr);
\draw(x2)[fill=white] circle(\vr);
\draw(x3)[fill=white] circle(\vr);
\draw(x4)[fill=white] circle(\vr);
\end{scope}

\begin{scope}[xshift=4cm, yshift=0cm] 
\coordinate(x1) at (1,3);
\coordinate(x2) at (2,4);
\coordinate(x3) at (1,5);
\coordinate(x4) at (0,4);
\draw (x1) -- (x2) -- (x3) -- (x4) --(x1);
\draw(x1)[fill=white] circle(\vr);
\draw(x2)[fill=white] circle(\vr);
\draw(x3)[fill=white] circle(\vr);
\draw(x4)[fill=white] circle(\vr);
\end{scope}

\begin{scope}[xshift=6cm, yshift=0cm] 
\coordinate(x1) at (1,3);
\coordinate(x2) at (2,4);
\coordinate(x3) at (1,5);
\coordinate(x4) at (0,4);
\draw (x1) -- (x2) -- (x3) -- (x4) --(x1);
\draw(x1)[fill=white] circle(\vr);
\draw(x2)[fill=white] circle(\vr);
\draw(x3)[fill=white] circle(\vr);
\draw(x4)[fill=white] circle(\vr);
\end{scope}

\begin{scope}[xshift=8cm, yshift=0cm] 
\coordinate(x1) at (1,3);
\coordinate(x2) at (2,4);
\coordinate(x3) at (1,5);
\coordinate(x4) at (0,4);
\draw (x1) -- (x2) -- (x3) -- (x4) --(x1);
\draw(x1)[fill=white] circle(\vr);
\draw(x2)[fill=white] circle(\vr);
\draw(x3)[fill=white] circle(\vr);
\draw(x4)[fill=white] circle(\vr);
\end{scope}

\begin{scope}[xshift=10cm, yshift=0cm] 
\coordinate(x1) at (1,3);
\coordinate(x2) at (2,4);
\coordinate(x3) at (1,5);
\coordinate(x4) at (0,4);
\coordinate(w) at (3,4);
\node at (2,4.3) {$v$};
\node at (3,4.3) {$u$};
\draw (x1) -- (x2) -- (x3) -- (x4) --(x1);
\draw (x2) -- (w);
\draw(x1)[fill=white] circle(\vr);
\draw(x2)[fill=white] circle(\vr);
\draw(x3)[fill=white] circle(\vr);
\draw(x4)[fill=white] circle(\vr);
\draw(w)[fill=white] circle(\vr);
\end{scope}

\begin{scope}[xshift=0cm, yshift=0cm] 
\coordinate(x1) at (0.5,0);
\coordinate(x2) at (1.5,0);
\coordinate(x3) at (2,1.1);
\coordinate(x4) at (1,1.8);
\coordinate(x5) at (0,1.1);
\draw (x1) -- (x2) -- (x3) -- (x4) --(x5) -- (x1);
\draw(x1)[fill=black] circle(\vr);
\draw(x2)[fill=white] circle(\vr);
\draw(x3)[fill=white] circle(\vr);
\draw(x4)[fill=white] circle(\vr);
\draw(x5)[fill=black] circle(\vr);
\end{scope}

\begin{scope}[xshift=2cm, yshift=0cm] 
\coordinate(x1) at (0.5,0);
\coordinate(x2) at (1.5,0);
\coordinate(x3) at (2,1.1);
\coordinate(x4) at (1,1.8);
\coordinate(x5) at (0,1.1);
\draw (x1) -- (x2) -- (x3) -- (x4) --(x5) -- (x1);
\draw(x1)[fill=white] circle(\vr);
\draw(x2)[fill=white] circle(\vr);
\draw(x3)[fill=white] circle(\vr);
\draw(x4)[fill=white] circle(\vr);
\draw(x5)[fill=white] circle(\vr);
\end{scope}

\begin{scope}[xshift=4cm, yshift=0cm] 
\coordinate(x1) at (0.5,0);
\coordinate(x2) at (1.5,0);
\coordinate(x3) at (2,1.1);
\coordinate(x4) at (1,1.8);
\coordinate(x5) at (0,1.1);
\draw (x1) -- (x2) -- (x3) -- (x4) --(x5) -- (x1);
\draw(x1)[fill=white] circle(\vr);
\draw(x2)[fill=white] circle(\vr);
\draw(x3)[fill=white] circle(\vr);
\draw(x4)[fill=white] circle(\vr);
\draw(x5)[fill=white] circle(\vr);
\end{scope}

\begin{scope}[xshift=6cm, yshift=0cm] 
\coordinate(x1) at (0.5,0);
\coordinate(x2) at (1.5,0);
\coordinate(x3) at (2,1.1);
\coordinate(x4) at (1,1.8);
\coordinate(x5) at (0,1.1);
\draw (x1) -- (x2) -- (x3) -- (x4) --(x5) -- (x1);
\draw(x1)[fill=white] circle(\vr);
\draw(x2)[fill=white] circle(\vr);
\draw(x3)[fill=white] circle(\vr);
\draw(x4)[fill=white] circle(\vr);
\draw(x5)[fill=white] circle(\vr);
\end{scope}

\begin{scope}[xshift=8cm, yshift=0cm] 
\coordinate(x1) at (0.5,0);
\coordinate(x2) at (1.5,0);
\coordinate(x3) at (2,1.1);
\coordinate(x4) at (1,1.8);
\coordinate(x5) at (0,1.1);
\draw (x1) -- (x2) -- (x3) -- (x4) --(x5) -- (x1);
\draw(x1)[fill=white] circle(\vr);
\draw(x2)[fill=white] circle(\vr);
\draw(x3)[fill=white] circle(\vr);
\draw(x4)[fill=white] circle(\vr);
\draw(x5)[fill=white] circle(\vr);
\end{scope}

\begin{scope}[xshift=10cm, yshift=0cm] 
\coordinate(x1) at (0.5,0);
\coordinate(x2) at (1.5,0);
\coordinate(x3) at (2,1.1);
\coordinate(x4) at (1,1.8);
\coordinate(x5) at (0,1.1);
\coordinate(w) at (3,1.1);
\draw (x1) -- (x2) -- (x3) -- (x4) --(x5) -- (x1);
\draw (x3) -- (w);
\draw(x1)[fill=white] circle(\vr);
\draw(x2)[fill=white] circle(\vr);
\draw(x3)[fill=white] circle(\vr);
\draw(x4)[fill=white] circle(\vr);
\draw(x5)[fill=white] circle(\vr);
\draw(w)[fill=black] circle(\vr);
\end{scope}

\end{tikzpicture}
\caption{Graphs $G_{6,4}$ and $G_{6,5}$ and their largest dual general position sets.}
	\label{fig:G64-G65}
\end{center}
\end{figure}

The converse of Proposition~\ref{prop:dual=1} is not true. For a sporadic example consider the graph $G$ obtained from $C_4$ by attaching a pendant vertex to one of the vertices of $C_4$. Then $\gpd(G) = 2 \ne 1$ but $s(G) = 1$. An infinite family of such examples is the following. For integers $k\geq 1$ and $\ell\geq 4$, let $G_{k,\ell}$ be the graph consisting of a chain of $k$ cycles $C_\ell$ which share vertices such that for an intermediate $C_\ell$ its vertices of degree $3$ are its diametral vertices. Finally attach a pendant vertex to a degree $2$ vertex of the last cycle. In Fig.~\ref{fig:G64-G65} the graphs $G_{6,4}$ and $G_{6,5}$ are shown. 

\begin{proposition}
\label{prop:dual=3}
If $k\geq 1$ and  $\ell\geq 4$, then
$$\gpd(G_{k,\ell}) = \left\{
\begin{array}{ll}
1; & \ell\geq 6\,,\\
2; & \ell=4\,,\\
3; & \ell=5\,.\\
\end{array}\right.
$$
\end{proposition}

\proof
Let $e = uv$ be the pendent edge of $G_{k,\ell}$, where $u$ is the pendant vertex. Since $u$ is a simplicial vertex of $G_{k,\ell}$, Corollary~\ref{cor:simplicial-are-dual} implies $\gpd(G_{k,\ell})\geq 1$.

Assume first $\ell\geq 6$. Since $g(G_{k,\ell})\geq 6$ and each edge of $G_{k,\ell}-u$ is \name, no vertex of $G_{k,\ell}-u$ lies in a dual general position set of $G_{k,\ell}$, cf.\ Corollary~\ref{cor:gpd=0 of girth 6}. Hence we have $\gpd(G_{k,\ell})\leq 1$ and thus $\gpd(G_{k,\ell}) = 1$.

Assume second $\ell = 4$. Let $i, i', j, j'$ be the vertices of the first $C_4$ as indicated in Fig.~\ref{fig:G64-G65}. Since $\{i,i'\}$ is a general position set of $G_{k,\ell}$ and $G_{k,\ell}-\{i,i'\}$ is convex, Theorem~\ref{thm:dual-characterization} implies that $\gpd(G_{k,\ell})\geq 2$. Suppose on the contrary that $\gpd(G_{k,\ell})\geq 3$ and let $X$ be an arbitrary $\gpd$-set of $G_{k,\ell}$. 

Assume that $X$ contains a vertex $x$ with $\deg_{G_{k,\ell}}(x) = 4$. Then at least three neighbors of $x$ lie in $X$ for, otherwise two neighbors of $x$ are not $X$-positionable. But then again two neighbors of $x$ are not $X$-positionable. By a parallel argument we also see that $X$ does not contain the vertex of degree $3$. It follows that each vertex of $X$ is of degree at most $2$. 

We claim that $u\notin X$. Since $|X|\geq 3$, there exist $y, y'\in X$ such that $\deg_{G_{k,\ell}}(y) = \deg_{G_{k,\ell}}(y') = 2$. (It is possible $y,y'\in\{i,i',j\}$.) If $y$ and $y'$ are adjacent, then we may without loss of generality assume that $y=i$ and $y' = i'$. But then $i'$ and $u$ are not $X$-positionable, so this cannot happen. Assume next that $y$ and $y'$ are not adjacent. Let $z$ and $z'$ be the two neighbors of $y$. It is clear that $z$ and $z'$ are not $X$-positionable if $y$ and $y'$ lie on the same cycle $C_4$. And if $y$ and $y'$ are not on the same $4$-cycle, then either $y$ lies on a shortest $u,y'$-path or $y'$ lies on the shortest $u,y$-path. We conclude that indeed $u\notin X$. 

We have thus proved that $\deg_{G_{k,\ell}}(x) = 2$ for each $x\in X$. Let $x$ and $x'$ be two arbitrary vertices from $X$. By the same argument as used in the previous paragraph we get that $xx'\in E(G_{k,\ell})$. As there are only two such edges possible, that is, $ii'$ and $i'j$, and as $\{i,i',j\}$ is not a dual general position set, we can conclude that $\{i,i'\}$ is a largest dual general position set and so $\gpd(G_{k,\ell}) = 2$.

The argument for the case $\ell = 5$ is similar and left to the reader. 
\qed

\subsection{Graphs $G$ with $\gpd(G) \ge 2$}
\label{sec:gpd=2}

We next consider when a set of cardinality two forms a dual general position set. In the next result we deduce a characterization (additional to the one of Theorem~\ref{thm:dual-characterization}) for two adjacent vertices.

\begin{theorem}
\label{thm:adjacent two vertices}
If $x$ and $y$ are two adjacent vertices of a graph $G$, then the following statements are equivalent.
\begin{enumerate}
\item[(i)]$\{x,y\}$ is a dual general position set of $G$;
\item[(ii)] $G-\{x,y\}$ is convex;
\item[(iii)] For each $u,v\in N_G(x)\cup N_G(y)$ we have $d_{G}(u,v)\leq 2$, and the graphs $G[N_G(x)-\{y\}]$ and $G[N_G(y)-\{x\}]$ are complete.
\end{enumerate}
\end{theorem}

\proof
Let $X=N_G(x)-\{y\}$ and $Y=N_G(y)-\{x\}$.

$(i) \Rightarrow (ii)$: If $\{x,y\}$ is a dual general position set of $G$, it follows from Theorem~\ref{thm:dual-characterization} that $G-\{x,y\}$ is convex.

$(ii) \Rightarrow (iii)$: Assume that $G-\{x,y\}$ is convex. Then $G[X]$ is complete. Indeed, for otherwise two nonadjacent vertices $x'$ and $x''$ from
$X$ are not $\{x,y\}$-positionable. Analogously, $G[Y]$ is complete. Furthermore, if $u\in X$ and $v\in Y$, then $d_G(u,v) \le 3$. But if $d_G(u,v) = 3$, then $G-\{x,y\}$ is not convex, hence we conclude that $d_{G}(u,v)\leq 2$ for any $u,v\in N_G(x)\cup N_G(y)$.

$(iii) \Rightarrow (i)$:
Let $1\leq d_{G}(u,v)\leq 2$ for any two vertices $u,v\in N_G(x)\cup N_G(y)$, and let $G[X]$ and $G[Y]$ be complete. Set $G'=G-\{x,y\}$. Consider two arbitrary distinct vertices $p$ and $q$ from $G'$. We claim that no shortest $p,q$-path passes $x$ or $y$. Let $Q$ be an arbitrary shortest $p,q$-path. Suppose first that $V(Q) \cap \{x,y\} = \{x\}$. Then $Q$ contains two neighbors of $x$ but since $G[X]$ is complete, this is a contradiction. The case when $V(Q) \cap \{x,y\} = \{y\}$ is ruled out analogously. In the remaining case suppose that $V(Q) \cap \{x,y\} = \{x,y\}$. Then $Q$ contains a subpath $x', x, y, y'$, where $x'\in X$ and $y'\in Y$. But by our assumption $d_G(x',y') \le 2$, a contradiction with the assumption that $Q$ is a shortest path. We conclude that $u$ and $v$ are $\{x,y\}$-positionable and consequently $\{x,y\}$ is a dual general position set.
\qed

A result parallel to Theorem~\ref{thm:adjacent two vertices} for two nonadjacent vertices is simpler and reads as follows.

\begin{proposition}
\label{prop:non-adjacent two vertices}
Let $x$ and $y$ be two non-adjacent vertices of a graph $G$. Then the set $\{x,y\}$ is a dual general position set if and only if $x$ and $y$ are simplicial vertices.
\end{proposition}

\proof
Assume first that $\{x,y\}$ is a dual general position set. Then $x$ is simplicial for otherwise there exist two neighbors $u$ and $v$ of $x$ such that $d_G(u,v) = 2$. By our assumption, $u\ne y$ and $v\ne y$, but then $u$ and $v$ are not $\{x,y\}$-positionable. Analogously $y$ is simplicial. The reverse implication follows by Corollary~\ref{cor:simplicial-are-dual}.
\qed

\section{The variety in Cartesian products}
\label{sec:in-Cartesian}

To determine the general position number of Cartesian product graphs is a difficult problem and has been already widely investigated. It took several intermediately steps before the general position number of integer lattices (alias Cartesian products of a finite number of paths) has been determined~\cite{Klavzar-Rus-2021}. Bounds on the general position number of Cartesian products of arbitrary graphs were independently proved in~\cite{Ghorbani-2021, Tian-2021}. Special Cartesian products were studied in ~\cite{Tian-2021b} (products of two trees), in~\cite{Klavzar-2021} (products of paths and cycles), and in~\cite{Klavzar-2021, Korze-Vesel-2023} (products of two cycles).

In this section we determine the total general position number, the outer general position number, and the dual general position number of arbitrary Cartesian products. For this purpose, we need some additional notation and terminology on Cartesian products. Let $G$ and $H$ be graphs and consider $G\cp H$. Given a vertex $h\in V(H)$, the subgraph of $G\cp H$ induced by the set of vertices $\{(g,h):\ g\in V(G)\}$, is a {\em $G$-layer} and is denoted by $G^h$. $H$-layers $^gH$ are defined analogously. Each $G$-layer and each $H$-layer is isomorphic to $G$ and $H$, respectively.
If $X\subseteq V(G\cp H)$, the {\em projection} $p_G(X)$ of $X$ to $G$ is the set $\{g\in V(G):\ (g,h)\in X\ {\rm for\ some}\ h\in V(H)\}$. The projection $p_H(X)$ of $X$ to $H$ is defined analogously.

For total general position sets, Corollary~\ref{cor:gpt=0} implies:

\begin{corollary}
If $G$ and $H$ are connected graphs of order at least $2$, then $\gpt(G\cp H) = 0$.
\end{corollary}

\proof
It is straightforward to see that  $G\cp H$ contains no simplicial vertices.
\qed

For outer general position sets, the following result is useful.

\begin{theorem} {\rm \cite[Theorem 3]{Rodriguez-2014}}
\label{thm:cartesian-direct}
If $G$ and $H$ are connected graphs, each of order at least $2$, then $(G\cp H)\sr \cong G\sr \times H\sr$.
\end{theorem}

\begin{theorem}
\label{thm:cartesian product for outer}
If $G$ and $H$ are two connected graphs, each of order at least $2$, then $$\gpo(G\cp H) = \min\{\gpo(G),\gpo(H)\}.$$
\end{theorem}

\proof
By Theorem~\ref{thm:characterize outer} we have $\gpo(G\cp H) = \omega ((G\cp H)\sr)$.
Theorem~\ref{thm:cartesian-direct} implies that $\omega ((G\cp H)\sr) = \omega (G\sr\times H\sr)$.
Since $\omega(G\times H) = \min \{\omega(G), \omega(H)\}$, see~\cite[Exercise 16.1]{HIK-2011}, we have $\gpo(G\cp H) =\min\{\omega (G\sr),\omega (H\sr)\}$.
Using Theorem~\ref{thm:characterize outer} again we conclude that
$\gpo(G\cp H) = \min\{\gpo(G),\gpo(H)\}$.
\qed

\begin{corollary}
\label{cor:complete graph}
If $m, n\geq 3$, then $\gpo(K_m\cp K_n) = \min \{m, n\}$.
\end{corollary}

For dual general position sets, the following result that characterizes convex subgraphs of Cartesian product graphs will be applied.

\begin{lemma} {\rm \cite[Lemma 6.5]{HIK-2011}}
\label{lem:subgraphs convex in cartesian product}
A subgraph $W$ of $G=G_1\cp\cdots\cp G_k$ is convex if and only if $W=U_1\cp\cdots\cp U_k$, where each $U_i$ is convex in $G_i$.
\end{lemma}

\begin{theorem}
\label{thm:cartesian product for dual}
Let $G$ and $H$ be two graphs, each of order at least $2$. Then $\gpd(G\cp H) > 0$ if and only if one factor is complete and the other factor has a simplicial vertex. Moreover, 
\begin{enumerate}
\item[(i)] $\gpd(K_n\cp K_m) = \max \{n,m\}$, and 
\item[(ii)] if $H$ is not complete and contains a simplicial vertex, then $\gpd(K_n\cp H) = n$.
\end{enumerate}
\end{theorem}

\proof
Set $V(G)=[n]$ and $V(H)=[m]$. So we have $V(G\cp H) = [n]\times [m]$.

Assume first that $\gpd(G\cp H) > 0$ and let $X$ be a dual general position set of $G\cp H$. Then $X$ is clearly a general position set of $G\cp H$ and hence by Theorem~\ref{thm:dual-characterization}, $K = (G\cp H) - X$ is a convex subgraph of $G\cp H$. Consequently, by Lemma~\ref{lem:subgraphs convex in cartesian product} we infer that $K = G' \cp H'$, where $G'$ is convex in $G$ and $H'$ is convex in $H$. Suppose that $G'$ and $H'$ are proper subgraphs of $G$ and $H$, respectively. Then there exist vertices $g\in V(G)\setminus V(G')$ and $h\in V(H)\setminus V(H')$, such that $g$ has a neighbor $g'\in V(G)$ and $h$ has a neighbor $h'\in V(H)$. Then $(g',h), (g,h), (g,h')$ is an induced path of $G\cp H$, where all its vertices are from $X$, a contradiction.

By the above contradiction, $p_G(X) = V(G)$ or  $p_H(X) = V(H)$. We may without loss of generality assume that $p_G(X) = V(G)$. Since layers in Cartesian products are convex, this implies that $V(G)$ is a general position set in $G$ which in turn implies that $G$ is a complete graph. Moreover, since $X$ is a general position set we also see that $|p_H(X)| = 1$ and let $h$ be the unique vertex of $H$ to which $X$ projects. Now, if $h$ is not a simplicial vertex of $H$, then there exist vertices $h', h''\in N_H(h)$ such that $h'h''\notin E(H)$. But then $K$ is clearly not convex, hence $h$ must be a simplicial vertex of $H$.

To complete the argument we observe that if $G$ is complete and $h\in V(H)$ is a simplicial vertex, then by Theorem~\ref{thm:dual-characterization} the set $V(G)\times \{h\}$ is a dual general position set.

The two formulas in the cases when $\gpd(G\cp H) > 0$ follow directly from the above discussion.
\qed

In~\cite[Theorem~3.8]{Tian-2021} it is proved that $\gp(K_m \cp K_n) = m + n - 2$. Combining this result with Corollary~\ref{cor:complete graph} and Theorem~\ref{thm:cartesian product for dual} we see that the four  general position invariants considered can vary arbitrary. For instance, for any $n\ge 3$ we have:
\begin{align*}
\gp(K_n\cp K_{2n}) & = 3n - 2,\\
\gpd(K_n\cp K_{2n}) & = 2n, \\
\gpo(K_n\cp K_{2n}) & =  n, \\
\gpt(K_n\cp K_{2n}) & = 0.
\end{align*}

\section{Concluding remarks}\label{sec:conclusion}

In this paper we have introduced the variety of general position sets. We have completely described the total general position sets and the outer general position sets. On the other hand, we have observed that the dual general position sets are not hereditary. This fact makes the investigation of dual general position sets quite tricky. For instance, we have given a sufficient condition for a graph $G$ to satisfy $\gpd(G) = 0$, yet we were not been able to characterize graphs $G$ with this property. The same problems remains open for the case $\gpd(G) = 1$. 

We have seem that if $G$ is a block graph, then $\gp(G) = \gpo(G) = \gpd(G) = \gpt(G)$. It would be an interesting project to characterize the graphs for which this holds. Moreover, the same question can be posed for each subsets of the involved invariants, for instance, to characterize the graphs $G$ for which $\gp(G) = \gpd(G)$ holds.

\section*{Acknowledgments}

Many thanks to James Tuite, who suggested to investigate the variety of the general position sets in parallel to the variety of mutual-visibility sets. This work was supported by the Slovenian Research Agency ARIS (research core funding P1-0297 and projects J1-2452, N1-0285). 

\section*{Declaration of interests}
 
The authors declare that they have no conflict of interest. 

\section*{Data availability}
 
Our manuscript has no associated data.


\end{document}